\documentclass[12pt]{amsart}
\usepackage[mathscr]{eucal}
\usepackage{amsmath,amssymb,hyperref}
\usepackage{graphics}

\newtheorem{thm}{Theorem}[section]

\newtheorem{cor}[thm]{Corollary}
\newtheorem{prop}[thm]{Proposition}

\newtheorem{theorem}{Theorem}
\newtheorem{corollary}[theorem]{Corollary}

\begin{document}

\title[A Preparation Theorem for foliations]{A Preparation Theorem 
for codimension one foliations}
\author[Frank Loray]{Frank LORAY}
\address{Frank LORAY (Charg\'e de Recherches au CNRS)
\hfill\break IRMAR, Campus de Beaulieu, 35042 Rennes Cedex (France)}
\email{frank.loray@univ-rennes1.fr}
\date{February 2004 (preliminary version)}
\dedicatory{to C\'esar Camacho for his $60^{\text{th}}$ birthday}
\thanks{}

\maketitle

\begin{abstract}After gluing foliated complex manifolds,
we derive a preparation-like theorem for singularities
of codimension one foliations and planar vector fields
(in the real or complex setting). 
Without computation, we retrieve and improve results
of Levinson-Moser for functions, Dufour-Zhitomirskii 
for non degenerate codimension $1$ foliations
(proving in turn the analyticity),
Strozyna-Zoladek for non degenerate planar vector fields
and Brjuno-\'Ecalle for saddle-node foliations in the plane.
\end{abstract}

\section*{Introduction}\label{S:intro}

We denote by $(\underline{z},w)$ the variable of $\mathbb{C}^{n+1}$,
$\underline{z}=(z_1,\ldots,z_n)$, for $n\ge1$. Recall that a germ 
of (non identically vanishing) holomorphic $1$-form 
$$
\Theta=f_1(\underline{z},w)dz_1+\cdots+f_n(\underline{z},w)dz_n+g(\underline{z},w)dw
$$
$f_1,\ldots,f_n,g\in\mathbb{C}\{\underline{z},w\}$,
defines
a codimension $1$ singular foliation $\mathcal F$ (regular outside the zero-set 
of $\Theta$) if, and only if, it satisfies Frobenius integrability condition
$\Theta\wedge d\Theta=0$. 
Maybe dividing coefficients of $\Theta$ by a common factor,
the zero-set of $\Theta$ has codimension $2$ 
and the foliation $\mathcal F$ extends
as a regular foliation outside this sharp singular set.

\eject

Our main result is the

\begin{theorem}\label{T:prepfol}
Let $\Theta$ and $\mathcal F$ as above and assume 
that $g(\underline{0},w)$
vanishes at the order $k\in\mathbb{N}^*$ at $0$.
Then, up to analytic change of the $w$-coordinate
$w:=\phi(\underline{z},w)$, the foliation $\mathcal F$ is 
also defined by a $1$-form
$$\widetilde\Theta=P_1(\underline{z},w)dz_1+\cdots+P_n(\underline{z},w)dz_n+Q(\underline{z},w)dw$$
for $w$-polynomials $P_1,\ldots,P_n,Q\in\mathbb{C}\{\underline{z}\}[w]$ of degree $\le k$,
$Q$ unitary.
\end{theorem}

In new coordinates given by Theorem \ref{T:prepfol}, 
the singular foliation $\mathcal F$ extends
analytically along some infinite cylinder 
$\{\vert \underline z\vert<r\}\times\overline{\mathbb C}$
(where $\overline{\mathbb C}=\mathbb C\cup\{\infty\}$ 
stands for the Riemann sphere).
To prove this Theorem, we just do the converse.
Given a germ of foliation, we force its endless analytic continuation 
in one direction
by constructing it in the simplest way, gluing foliated manifolds into a 
foliated $\overline{\mathbb C}$-bundle. This is done in Section \ref{S:fol}.
The huge degree of freedom encountered during our construction 
can be used to preserve additional structure equipping the foliation. 
For instance, starting with the complexification of a real analytic foliation,
our gluing construction can be carried out preserving
the anti-holomorphic involution 
$(\underline z,w)\mapsto(\overline{\underline z},\overline w)$
so that our statement agree with the real setting.
By the same way, if one starts with a closed meromorphic $1$-form $\Theta$,
one can arrange so that $\Theta$ extends meromorphically
as well along the infinite cylinder (see Section \ref{S:closedforms}) 
and becomes itself rational in $w$.
In particular, in the case $\Theta=df$ is exact, 
we derive the following alternate Preparation Theorem.

\begin{theorem}[Levinson]\label{T:prepfunct}
Let $f(\underline{z},w)$ be a germ of holomorphic function at $(\underline 0,0)$ in $\mathbb{C}^{n+1}$
and assume that $f(\underline{0},w)$ vanishes at the order $k\in\Bbb N^*$ at $w=0$.
Then, up to an analytic change of coordinates,
the function germ $f$ becomes a unitary $w$-polynomial of degree $k$
$$f(\underline{z},w)=w^k+f_{k-1}(\underline{z}) w^{k-1}+\cdots+f_0(\underline{z})$$
where $f_0,\ldots,f_{k-1}\in\mathbb{C}\{\underline{z}\}$.
\end{theorem}

The difference with Weierstrass Preparation Theorem lies in the fact
that the usual invertible factor term (in variables $(\underline{z},w)$) 
is normalized to $1$ here;
the counter part is that a change of coordinates is needed.
This result was previously obtained by N. Levinson in \cite{Le} 
after an iterative procedure and reproved by J. Moser in \cite{Mo}
as an example illustrating KAM fast convergence. Similarly, we obtain that any germ of meromorphic
function is conjugated to a quotient of Weierstrass $w$-polynomials
(see Theorem \ref{T:merofunct}).

\eject

For $k=1$, Theorem \ref{T:prepfol} specifies as follows.

\begin{corollary}\label{T:fol}
Let $\Theta$ and $\mathcal F$ as in Theorem \ref{T:prepfol}
and assume that the linear part of $\Theta$ is not tangent
to the radial vector field.
Then, there exist local analytic coordinates $(\underline{z},w)$
in which the foliation $\mathcal F$ is defined by
$$\widetilde\Theta=df_0+wdf_1+wdw$$
where $f_0,f_1\in\mathbb{C}\{\underline z\}$ satisfy $df_0\wedge df_1=0$.
\end{corollary}

Following \cite{MaMo}, the functions $f_i$ factorize into a primitive function $f$
and the foliation $\mathcal F$ is actually the lifting 
of a foliation in the plane by the holomorphic map 
$\Phi:(\mathbb C^{n+1},\underline 0)\to(\mathbb C^2,0);
(\underline z,w)\mapsto(f(\underline z),w)$.
This normal form has been obtained in \cite{DuZh} 
by J.-P. Dufour and M. Zhitomirski 
after a formal change of coordinates
but the convergence was not proved.

In Theorem \ref{T:prepfol}, the $\overline{\mathbb C}$-fibration is constructed
simultaneously with the extension of the foliation $\mathcal F$
by gluing bifoliated manifolds. 
In dimension $2$, when $\mathcal F$ is defined by a vector field $X$,
it is still possible to extend $X$ on a $2$-dimensional tubular neighborhood $M$
of an embedded sphere $\overline{\mathbb C}$ but it is no more possible 
to construct the $\overline{\mathbb C}$-fibration at the same time.
Here, we need the Rigidity Theorem of V. I. Savelev (see \cite{Sa}):
{\it the germ of $2$-dimensional neighborhood of an embedded sphere having zero 
self-intersection is a trivial $\overline{\mathbb C}$-bundle over the disc}.
In Section \ref{S:vectfields}, we derive, 
for non degenerate singularities of vector fields

\begin{theorem}\label{T:vectfields}Let $X$ be a germ of analytic vector field
vanishing at the origin of $\mathbb{R}^2$ 
(resp. of $\mathbb{C}^2$).
Assume that its linear part is not radial.
Then, there exist local analytic coordinates $(x,y)$
in which 
$$X=(y+f(x))\partial_x+g(x)\partial_y$$
where $f,g\in\mathbb{R}\{x\}$ (resp. $f,g\in\mathbb{C}\{x\}$) vanish at $0$.
\end{theorem}

Denote by $\lambda_1,\lambda_2\in\mathbb C$ the eigenvalues of the vector field $X$:
we have $\lambda_1+\lambda_2=f'(0)$ and $\lambda_1\cdot\lambda_2=-g'(0)$.
In the case $\lambda_2=-\lambda_1$ (including nilpotent case $\lambda_i=0$), 
Theorem \ref{T:vectfields}
was obtained by E. Strozyna and H. Zoladek in \cite{StZo}.
They did prove the convergence
of an explicit iterative reduction process
after long and technical estimates.
In the case $\lambda_2/\lambda_1\not\in\mathbb R^-$,
Theorem \ref{T:vectfields} becomes just useless since
H. Poincar\'e and H. Dulac gave a unique and very simple polynomial 
normal form. 
In the remaining case,
taking in account the invariant curve of the vector field $X$,
we can specify our normal form as follows (see Section \ref{S:vectfields}
for a statement including nilpotent singularities).

\begin{corollary}\label{C:vectfields}Let $X$ be a germ of 
analytic vector field in the real or complex plane with eigenratio
$\lambda_2/\lambda_1\in\mathbb R^-$.
Then, there exist local analytic coordinates
in which the vector field $X$ takes the form
\begin{enumerate}
\item in the saddle case $\lambda_2/\lambda_1\in\mathbb{R}^-_*$
(with $\lambda_1,\lambda_2\in\mathbb R$ in the real case)
$$X=f(x+y)\left\{(\lambda_1 x\partial_x+\lambda_2 y\partial_y)+g(x+y)(x\partial_x+y\partial_y)\right\}$$
\item in the saddle-node case, say $\lambda_2=0$, $\lambda_1\not=0$
$$X=f(x)\left\{(\lambda_1 x+y)\partial_x+g(x)y\partial_y\right\}$$
\item in the real center case $\lambda_2=-\lambda_1=i\lambda$, $\lambda\in\mathbb R$
$$X=f(x)\left\{(-\lambda y\partial_x+\lambda x\partial_y) + g(x)(x\partial_x+y\partial_y)\right\}$$
\end{enumerate}
In each case, $f(0)=1$ and $g(0)=0$.
\end{corollary}

Orbital normal form (i.e. normal form for the induced foliation)
can be immediately derived just by setting $f\equiv 1$: coefficient $g$
stands for the moduli of the foliation.
Normal form (3) was also derived in \cite{StZo}.

In case (1), A. D. Brjuno proved in \cite{Br1} that the vector field $X$ 
is actually analytically linearisable for generic eigenratio 
$\lambda_2/\lambda_1\in\mathbb R^-$ 
(in the sense of the Lebesgue measure).
In this case, normal form (1) of 
Theorem \ref{C:vectfields} becomes just useless.
For the remaining exceptional values,
the respective works of J.-C. Yoccoz in the diophantine case
(see \cite{Yo} and \cite{PMYo}) and J. Martinet with J.-P. Ramis in the
resonant case $\lambda_2/\lambda_1\in\mathbb Q^-$
(see \cite{MaRa}) derive a huge moduli space for the analytic classification
of the induced foliations.
This suggests that most of the vector fields having such eigenvalues 
are not polynomial in any analytic coordinates. 
Moreover, at least in the resonant case, 
the analytic classification 
of all vector fields inducing
a given foliation gives rise to functional moduli as well 
(see \cite{GrVo}, \cite{MeVo} and \cite{Te}).
Thus, the functional parameters $f$ and $g$ appearing in our normal form
seem necessary in many cases.

Finally, one can shortly derive from (2) the {\it versal deformation}
$$X_f=x\partial_x+y^2\partial_y+yf(x)\partial_x,\ \ \ f\in\mathbb C\{x\},$$
of the saddle-node foliation $\mathcal F_0$ defined by $X_0=x\partial_x+y^2\partial_y$ 
(see \cite{versal}).
In other words, any germ of analytic deformation of $X_0$ without bifurcation
of the saddle-node point factorizes into the family above after
analytic change of coordinates and renormalization. Moreover, the derivative 
of Martinet-Ramis' moduli map 
at $X_0$ (see \cite{El}) is bijective.
When $f(0)=0$, one can even show that the form above is unique. 
This result was announced by A.D. Brjuno in \cite{Br2} 
and proved by J. \'Ecalle at the end of \cite{Ec}
using {\it mould theory} in the particular case $f'(0)=0$. 
We will detail it in a forthcoming paper \cite{versal}.

\eject

\section{Preparation Theorem for codimension $1$ foliations}\label{S:fol}

We first prove Theorem \ref{T:prepfol}. 
Let $\mathcal F_0$ denote the germ of singular foliation defined by an
integrable holomorphic $1$-form at $(\underline 0,0)\in\mathbb C^{n+1}$
$$\Theta_0=f_1(\underline{z},w)dz_1+\cdots+f_n(\underline{z},w)dz_n+g(\underline{z},w)dw,\ \ \
\Theta_0\wedge d\Theta_0=0$$
$f_1,\ldots,f_n,g\in\mathbb{C}\{\underline{z},w\}$ 
and assume $g(\underline 0,w)\not\equiv0$.
In particular, for $r>0$ small enough, the 
foliation $\mathcal F_0$ is well-defined on the vertical disc 
$\Delta_0=\{\underline 0\}\times\{\vert w\vert<r\}$, regular
and transversal to $\Delta_0$ outside $w=0$.

Consider in $\mathbb C^n\times\overline{\mathbb C}$
the vertical line $L=\{\underline 0\}\times\overline{\mathbb C}$ 
together with the covering given by $\Delta_0$ and another disc, say 
$\Delta_\infty=\{\underline 0\}\times\{\vert w\vert>r/2\}$.
Denote by $C=\Delta_0\cap\Delta_\infty$ the intersection corona.
By the flow-box Theorem, there exists a unique germ of diffeomorphism 
of the form
$$\Phi:(\mathbb C^{n+1},C)\to (\mathbb C^{n+1},C)\ ;\ 
 (\underline z,w)\mapsto (\underline z,\phi (\underline z,w)),\ \ \ 
 \phi (\underline 0,w)=w$$
conjugating $\mathcal F_0$ to the horizontal foliation $\mathcal F_\infty$
(defined by $\Theta_\infty=dw$) at the neighborhood of the corona $C$.
Therefore, after gluing the germs of complex manifolds
$(\mathbb C^n\times\overline{\mathbb C},\Delta_0)$ and
$(\mathbb C^n\times\overline{\mathbb C},\Delta_\infty)$
along the corona by means of $\Phi$, we obtain a germ of smooth
complex manifold $M$, $\dim(M)=n+1$, along a rational curve $L$ equipped with 
a singular holomorphic foliation $\mathcal F$. 
Moreover, the coordinate $\underline z$, which
is invariant under the gluing map $\Phi$, defines a germ of rational
fibration $\underline z:(M,L)\to(\mathbb C^n,\underline 0)$.
Following \cite{FiGr}, there exists a germ of submersion 
$w:(M,L)\to L\simeq\overline{\mathbb C}$ 
completing $\underline z$ into a system of trivializing coordinates
$(\underline z,w):(M,L)\to(\mathbb C^n,\underline 0)\times\overline{\mathbb C}$.
This system is unique up to permissible change
$$(\tilde{\underline z},\tilde w)=
\left(\phi(\underline z),{a(\underline z)w+b(\underline z)
\over c(\underline z)w+d(\underline z)}\right)$$
where $a,b,c,d\in\mathbb C\{\underline z\}$, $ad-bc\not\equiv0$, and 
$\phi\in\text{Diff}(\mathbb C^n,\underline 0)$.

At the neighborhood of any point $p\in L$, the foliation $\mathcal F$
is defined by a (non unique) germ of holomorphic $1$-form (respectively $\Theta_0$ 
or $\Theta_\infty$) or equivalently by a unique germ of meromorphic
$1$-form 
$$\Theta=R_1(\underline{z},w)dz_1+\cdots+R_n(\underline{z},w)dz_n+dw$$
where $R_i$ are meromorphic at $p$. By unicity, this $1$-form is actually
globally defined, on the neighborhood of $L$, and is therefore rational
in the variable $w$, i.e. all coefficients $R_i$ are quotients
of Weierstrass polynomials. 

\eject

\begin{figure}[htbp]
\begin{center}

\input{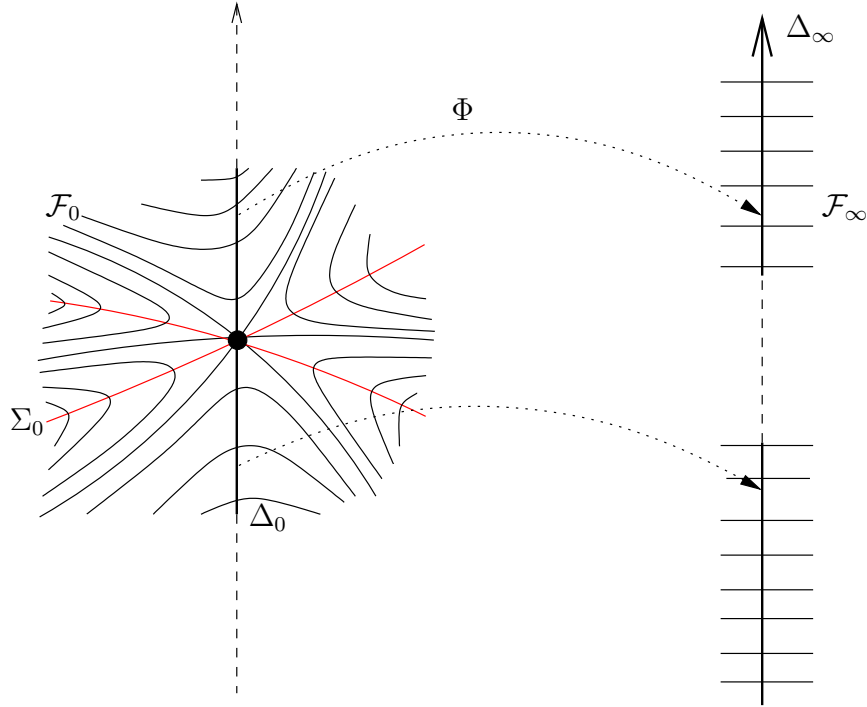}
 
\caption{Glueing construction}
\label{figure:1}
\end{center}
\end{figure}

Choose trivializing coordinates $(\underline z,w)$ so that
the singular point of $\mathcal F$ is still located at $w=0$.
The poles of $\Theta$ correspond
to tangencies between the foliation $\mathcal F$ and the rational fibration 
(counted with multiplicity). Denote by $\Sigma$ this divisor.
Since $\mathcal F_\infty$ is transversal to the rational fibration,
those poles come from the first chart, namely from the corresponding
tangency divisor
$$\Sigma_0=\{g(\underline z,w)=0\}$$
By assumption,
the total number of tangencies between $\mathcal F$ (or $\mathcal F_0$) 
and a fibre (close to $L$) is $k$. It follows that the
$w$-rational coefficients $R_i$ have exactly $k$ poles 
(counted with multiplicity) in restriction to each fiber. 
Therefore, if $Q$ denotes the unitary $w$-polynomial of degree $k$
defining $\Sigma$ and if one write $R_i={P_i\over Q}$ 
for $w$-polynomials $P_i$, the transversality of $\mathcal F$ 
with the fibration at $\{w=\infty\}$ implies that the $P_i$'s have at most
degree $k+2$ in variable $w$. Equivalently, $\mathcal F$ is defined by
$$\widetilde{\Theta}=\theta_{0}+\theta_1 w+\cdots+\theta_{k+2}w^{k+2}+Q(\underline{z},w)dw$$
for evident $1$-forms $\theta_0,\theta_1,\ldots,\theta_k$ on
$(\mathbb{C}^n,\underline 0)$ (depending only on $\underline z$).

\eject

After a permissible change of the $w$-coordinate,
one may assume that the line $\{w=\infty\}$ at infinity
is a leaf of the foliation (just straighten one leaf), 
i.e. that $\theta_{k+2}=0$. In fact, one may furthermore 
assume that the contact between $\mathcal F$ 
and the horizontal fibration $\{w=\text{constant}\}$ 
along the line $\{w=\infty\}$ has multiplicity $2$
(no linear holonomy along this leaf in the $w$-coordinate).
Indeed, make the change of coordinate $\tilde w=e^{-\int\theta_{k+1}}w$
($\theta_{k+1}$ is closed by integrability condition 
$\widetilde{\Theta}\wedge d\widetilde{\Theta}=0$).
In new coordinates, $\theta_{k+1}=0$ and Theorem \ref{T:prepfol}
is proved. Notice that we can further simplify the form $\widetilde{\Theta}$
by using the remaining possible changes of coordinates 
$\tilde{\underline z}=\phi(\underline z)$
and $\tilde w=w+b(\underline z)$.

\begin{figure}[htbp]
\begin{center}

\input{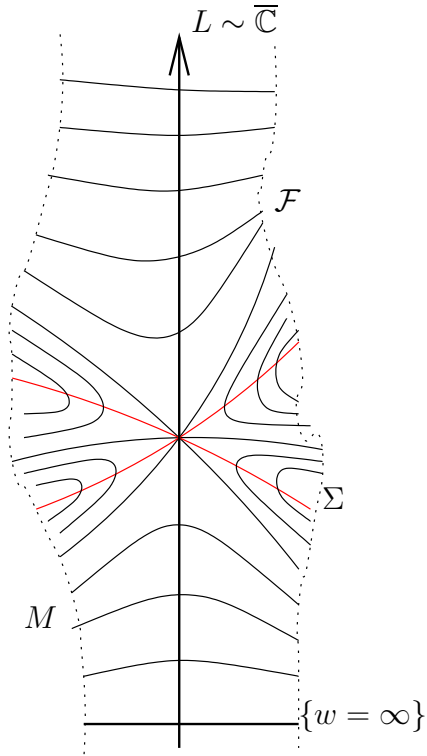}
 
\caption{Uniformisation}
\label{figure:3}
\end{center}
\end{figure}

\eject

We now prove Corollary \ref{T:fol}. Going back to the begining 
of the proof above, if the linear part of $\Theta_0$ is not 
tangent to the radial vector field,
up to a linear change of coordinates, one may assume that the tangency
set $\Sigma_0=\{g(\overline z,w)=0\}$ between the foliation $\mathcal F_0$
and the vertical fibration $\{\underline z=\text{constant}\}$
is smooth and transversal to the fibration. We are in the assumption
of Theorem \ref{T:prepfol} with $k=1$: up to a change of the $w$-coordinate,
one may assume that $\mathcal F$ is defined by
$\widetilde{\Theta}=\theta_0+w\theta_1+(w+f(\underline z))dw$
where $\theta_0$ and $\theta_1$ are holomorphic $1$-form depending only
on the $\underline z$-variable and $f\in\mathbb{C}\{\underline{z}\}$.
After translation $w:=w+f(\underline z)$ (notice that $f(\underline 0)=0$),
one may assume furthermore $f\equiv0$ and integrability condition
$\widetilde{\Theta}\wedge d\widetilde{\Theta}=0$ yields
$$\theta_0\wedge\theta_1=0,\ \ \ d\theta_0=0\ \ \ \text{and}\ \ \ d\theta_1=0.$$
After integration, we obtain $\theta_i=df_i$
for functions $f_i\in\mathbb{C}\{\underline{z}\}$ with tangency condition
$df_0\wedge df_1=0$ and Corollary \ref{T:fol} is proved. 
By \cite{MaMo}, there exists a primitive function
$f\in\mathbb{C}\{\underline{z}\}$ (with connected fibres)
through which $f_0$ and $f_1$ factorize: $f_i=\tilde f_i\circ f$ with
$\tilde f_i\in\mathbb C\{z\}$, $z$ a single variable.
Notice that we can further simplify the form $\widetilde{\Theta}$
by using the remaining possible changes of coordinate 
$\tilde{\underline z}=\phi(\underline z)$.

If we start with a real analytic foliation $\mathcal F_0$,
then its complexification is invariant under
the anti-holomorphic involution 
$(\underline z,w)\mapsto (\overline{\underline z},\overline w)$.
This involution obviously commutes with $\mathcal F_{\infty}$
and with the gluing map $\Phi$, defining by this way a germ 
of anti-holomorphic involution
$\Psi:(M,L)\to(M,L)$ on the resulting manifold preserving $\mathcal F$. 
In restriction to the coordinate $\overline z$, which is invariant under $\Phi$ 
and well defined on $M$, $\Psi$ induces the standart involution 
$\underline z\mapsto\overline{\underline z}$. Therefore, 
$\Psi(\underline z,w)=(\overline{\underline z},\psi(\underline z,w))$ 
where $\psi(\underline z,w)$ is, for fixed $\underline z$, a reflection 
with respect to a real circle. After an holomorphic change of $w$-coordinate,
$\psi(\underline z,w)=\overline w$ and the constructed foliation $\mathcal F$
is actually invariant by the standart involution. 
The unique meromorphic $1$-form defining $\mathcal F$ 
$$\Theta=R_1(\underline{z},w)dz_1+\cdots+R_n(\underline{z},w)dz_n+dw$$
satisfies
$\Psi^*\Theta=\overline{\Theta}$ and its coefficients
are actually real: $R_i\in\mathbb R\{\underline z\}(w)$.
This real form is obtained up to global change of coordinates commuting with
the standart involution, that is
$$(\tilde{\underline z},\tilde w)=
\left(\phi(\underline z),{a(\underline z)w+b(\underline z)
\over c(\underline z)w+d(\underline z)}\right)$$
where $a,b,c,d\in\mathbb R\{\underline z\}$, $ad-bc\not\equiv0$, and 
$\phi\in\text{Diff}(\mathbb R^n,\underline 0)$.

\eject

\section{Preparation Theorem for closed meromorphic $1$-forms}\label{S:closedforms}

For simplicity, we start with the case of (meromorphic) functions:

\begin{thm}\label{T:merofunct}
Let $f$ be a germ of meromorphic function at $(\underline 0,0)$ in $\mathbb{C}^{n+1}$
and assume that $f(\underline{0},w)$ is a well-defined and non constant germ of meromorphic
function. Then, up to analytic change of the $w$-coordinate $w:=\phi(\underline{z},w)$,
the function $f$ becomes a $w$-rational function
$$f(\underline{z},w)={f_0(\underline{z})+f_1(\underline{z})w+\cdots
+f_{k_0-1}(\underline{z})w^{k_0-1}+w^{k_0}\over
g_0(\underline{z})+g_1(\underline{z}) w+\cdots
+g_{k_\infty-1}(\underline{z})w^{k_\infty-1}+w^{k_\infty}}$$
where $k_0,k_\infty\in\mathbb N$ and $f_i,g_j\in\mathbb{C}\{\underline{z}\}$.
\end{thm}

\begin{proof}
Denote by $f_0(\underline z,w)$ the germ of meromorphic function
above and make a preliminary change of coordinate $\tilde w:=\varphi(w)$
such that $f_0(\underline{0},w)=w^l$, $l\in\mathbb Z^*$, or $1+w^l$, $l\in\mathbb N^*$. 
Then, proceed with the
underlying foliation $\mathcal F_0$ (defined by $f_0=\text{constant}$)
as in the proof of Theorem \ref{T:prepfol} in section \ref{S:fol}.
By construction, the function $f_0$ will glue automatically with the 
respective function
$f_\infty(\underline{z},w)=w^l$ or $1+w^l$ defining $\mathcal F_\infty$.
Therefore, the global foliation $\mathcal F$ is actually defined 
by a global meromorphic function $f$ on $M$. Again, $f$ is a quotient
of Weierstrass polynomials. In the case $f_0(\underline{0},w)=w^l$,
choose the $w$-coordinate such that the zero or pole of 
$f_\infty(\underline{z},w)=w^l$ still coincides with $\{w=\infty\}$.
Therefore, $k_0$ and $k_\infty$ respectively coincide with the number of 
zeroes and poles of $f_0$ restricted to a generic vertical line (close to $L$).
In the other case $f_0(\underline{0},w)=1+w^l$, we add $l$ simple zeroes 
in the finite part and a pole of order $l$ that can be straightened to $\{w=\infty\}$
as before. In this latter case, $l=k_0-k_\infty>0$ and $k_\infty$ is the number
of (zeroes or) poles of $f_0(\underline z,w)$ restricted to a generic vertical line. 
In any case, the leading terms $f_{k_0}$ and $g_{k_\infty}$ are non vanishing 
at $\underline z=\underline 0$ and can be normalized to $1$ by division
and a further change of coordinate $\tilde w=a(\underline z)w$. 
\end{proof}

The proof of Theorem \ref{T:prepfunct} immediately follows from setting 
$k=k_0>0$ and $k_\infty=0$ in the proof above.

\begin{prop}\label{L:closedforms}
Let $\Theta$ be a germ of closed meromorphic $1$-form at $(\underline 0,0)\in\mathbb{C}^{n+1}$
and assume that the vertical line $\{\underline z=\underline 0\}$ is not invariant
by the induced foliation. Then, up to analytic change of the $w$-coordinate $w:=\phi(\underline{z},w)$,
the closed form $\Theta$ takes the form
$$\Theta={P_1(\underline{z},w)dz_1+\cdots+P_n(\underline{z},w)dz_n+P(\underline{z},w)dw
\over Q(\underline{z},w)}$$
for $w$-polynomials $P,Q,P_1,\ldots,P_n\in\mathbb C\{\underline{z}\}[w]$.
\end{prop}

\begin{proof}By a preliminary change of the $w$-coordinate, one can
normalize the restriction of $\Theta$ to the vertical line into
one of the models
\begin{eqnarray*}
\Theta\vert_{L}&=&w^kdw\ \text{if}\ k\ge0\\
\Theta \vert_{L}&=&\lambda {dw\over w}\ \text{if}\ k=-1\\
\Theta \vert_{L}&=&\lambda {dw\over w^k(1-w)}\ \text{if}\ k<-1
\end{eqnarray*}
where $k\in\mathbb Z$ stands for the order of $\Theta\vert_{L}$ at $w=0$
and $\lambda\in\mathbb C$ denote the residue when $k\le-1$.
Then, defining the horizontal foliation $\mathcal F_\infty$
by the corresponding model $\Theta_\infty$ above 
(viewed as a $1$-form in variables $(\underline z,w)$),
we proceed gluing the foliations and the $1$-forms as we did with
functions in the previous proof. 
If $k_0$ and $k_\infty$ denote the respective number of zeroes and poles
of $\Theta_0$ in restriction to a generic vertical line, 
then the numerator and denominator have respective degree $k_0$ and 
$k_\infty$ if $k_0-k_\infty\ge-1$ and
$k_0$ and $k_\infty+1$ if $k_0-k_\infty<-1$.
\end{proof}

\section{Non degenerate vector fields in the plane}\label{S:vectfields}

We prove Theorem \ref{T:vectfields} and deduce Corollary \ref{C:vectfields}.
Let $X_0$ be a germ of analytic vector field at $(0,0)\in\mathbb C^2$
$$X_0=f(z,w)\partial_z+g(z,w)\partial_w$$
vanishing at $(0,0)$ with a non radial linear part
$$\text{lin}(X_0)=(az+bw)\partial_z+(cz+dw)\partial_w
=\left(\begin{matrix}a&b\\c&d\end{matrix}\right)
\not=\left(\begin{matrix}\lambda&0\\0&\lambda\end{matrix}\right)$$
(in particular, it is assumed that the linear part is not the zero matrix).
One can find linear coordinates in which
$$\text{lin}(X_0)=\left(\begin{array}{cc}
  0 &  1 \\
  \alpha &  \beta \\
\end{array}\right)+\cdots$$
where $-\alpha$ and $\beta$ respectively stand for the product and the sum of the eigenvalues
$\lambda_1$ and $\lambda_2$. The eigenvector corresponding to $\lambda_i$ is $(1,\lambda_i)$;
in the case $\lambda_1=\lambda_2$, we note that the matrix above is not diagonal.
After a change of the $w$-coordinate of the form $w:=\varphi(w)$, we may assume that
restriction of $f(z,w)$ to the vertical line $\{z=0\}$ takes the form $f(0,w)=w$.
Similarly to the proof of Theorem \ref{T:prepfol} in section \ref{S:fol},
we consider in $\mathbb C\times\overline{\mathbb C}$
the vertical line $L=\{0\}\times\overline{\mathbb C}$ 
together with the covering given by 
$$\Delta_0=\{0\}\times\{\vert w\vert<r\}\ \ \ \text{and}\ \ \ 
\Delta_\infty=\{0\}\times\{\vert w\vert>r/2\}$$
and we denote by $C=\Delta_0\cap\Delta_\infty$ the intersection corona.

\eject

If $r>0$ is small enough, the vector fields $X_0$ is well defined on the neighborhood
of the closed disc $\overline{\Delta_0}$ and transversal to it outside $w=0$.
By Rectification Theorem, there exists a unique germ of diffeomorphism 
of the form
$$\Phi:(\mathbb C^2,C)\to (\mathbb C^2,C),\ \ \ 
 \Phi (0,w)=(0,w)$$
conjugating $X_0$ to the horizontal vector field $X_\infty=w\partial_z$.
After gluing the germs of complex surfaces
$(\mathbb C\times\overline{\mathbb C},\Delta_0)$ and
$(\mathbb C\times\overline{\mathbb C},\Delta_\infty)$
along the corona by means of $\Phi$, we obtain a germ of smooth
complex surface $M$ along a rational curve $L$ equipped with 
a meromorphic vector field $X$. 
Since the $\partial_z$-component of $X_0$ agree with $w\partial_z$ along
$L$, it follows that the Jacobian of the gluing map $\Phi$ takes the form
$$D_{(0,w)}\Phi=\left(\begin{array}{cc}
  1 &  0 \\
 * &  1 \\
\end{array}\right)$$
and the embedded rational curve $L$ has zero self-intersection.
Following \cite{Sa}, $L$ is the regular fiber of a germ of a trivial fibration
on $M$, i.e. there exist global coordinates 
$(z,w):(M,L)\to(\mathbb C,\underline 0)\times\overline{\mathbb C}$
sending $L$ onto $\{z=0\}$. The vector field $X$ has exactly one isolated zero, 
say $(z,w)=(0,0)$,
and a simple pole along a trajectory (given by $X_\infty=w\partial_z$
in the second chart) that we may assume
still given by $\{w=\infty\}$. The tangency divisor $\Sigma$ between
the induced foliation $\mathcal F$ and the fibration $\{z=\text{constant}\}$
still is a smooth curve intersecting the fiber $\{z=0\}$ at the singular point
$(z,w)=(0,0)$ without multiplicity. Indeed, the Jacobian of the change of coordinates
(from the first chart to the global coordinates) at the singular point is fixing
the $w$-direction, so that the linear part of the vector field takes the form
$$X=\left(\begin{array}{cc}
  a &  b \\
 c &  d \\
\end{array}\right)+\cdots,\ \ \ b\not=0.$$
As in the proof of Theorem \ref{T:prepfol}, one may choose the (global) $w$-coordinate
so that the foliation has a contact of order $2$ with the horizontal
foliation $\{w=\text{constant}\}$ along the polar trajectory $\{w=\infty\}$
and the tangency set $\Sigma=\{w=0\}$ is horizontal as well.
Therefore, the vector field $X$ writes
$$X=f(z)w\partial_z+(g_0(z)+wg_1(z))\partial_w$$
for germs $f,g_0,g_1\in\mathbb C\{z\}$. Indeed,
the coefficients of $X=P(z,w)\partial_z+Q(z,w)\partial_w$
become automatically rational in the $w$-variable.
Since the unique pole of $X$ is simple and located at $\{w=\infty\}$,
$P$ and $Q$ are in fact polynomials of maximal degree $1$ and $3$
(notice that $\partial_w$ has a double zero at $\{w=\infty\}$).
Finally, conditions on tangency and polar sets imply the special form above.

\eject

\begin{figure}[htbp]
\begin{center}

\input{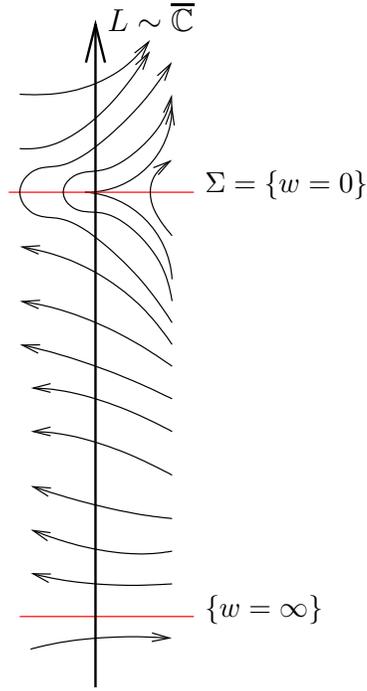}
 
\caption{Normalisation}
\label{figure:2}
\end{center}
\end{figure}

By a change of $z$-coordinate, we may furthermore assume $f(z)\equiv1$
($f(0)=1\not=0$). Automatically, the linear part of $X$ in the new coordinates writes
$$X=\left(\begin{array}{cc}
  0 &  1 \\
  \alpha &  \beta \\
\end{array}\right)+\cdots$$
i.e. $g_0(z)=\alpha z+\cdots$ and $g_1(z)=\beta+\cdots$ 
where dots mean higher order terms.
Finally, the form $X=(w+f(z))\partial_z+g(z)\partial_w$
is derived after the last change of coordinate $\tilde w:=w-f$ 
where $f'(z)=g_1(z)$, $f(0)=0$.

If we start with a real vector field $X_0$, 
then the anti-holomorphic involution 
$(z,w)\mapsto (\overline{z},\overline w)$
commutes with the gluing map $\Phi$ 
(mind that $X_\infty=w\partial_z$ is also real)
and induces a germ of anti-holomorphic involution
$\Psi:(M,L)\to(M,L)$ on the resulting surface 
satisfying $\overline{\Psi^*X}=X$. 
By Blanchard's Lemma, $\Psi$ preserves the rational fibration.
In restriction to the coordinate $z$, $\Psi$ is a regular
anti-holomorphic involution and is obviously holomorphically
conjugated to the standart one $z\mapsto\overline{z}$. 
Finally, after holomorphic change of $w$-coordinate, 
$\Psi(z,w)=(\overline z,\overline w)$ and $X$ has real
coefficients.

\eject

\begin{cor}\label{C:vectfieldsbis}Let $X$ be a germ of
analytic vector field as in Theorem \ref{T:vectfields}.
Then, by a further change of (complex or real) analytic coordinates,
one may assume that we are in one of the following cases
\begin{enumerate}
\item $X$ has an invariant curve of the form $C:\{w^2-z^k=0\}$ and
$$X=f(z)(2w\partial_z+k z^{k-1}\partial_w) + g(z)z^l(2z\partial_z+k w\partial_w),\ \ \ l+1\ge{k\over 2}\ge1$$
\item $X$ has an invariant curve of the form $C:\{w=0\}$ and
$$X=f(z)(w+z^k)\partial_z + g(z)z^lw\partial_w,\ \ \ l+1\ge k\ge1$$
\item $X$ is a real center or focus 
and
$$X=f(z)(-w\partial_z+k z^{2k-1}\partial_w) + g(z)z^l(z\partial_z+k w\partial_w),\ \ \ l+1\ge k\ge1$$
\end{enumerate}
where, in every case, $f(0)\not=0$.
\end{cor}

Saddles and saddle-nodes respectively correspond to cases 1 and 2.
For a complete discussion on the possible invariant curve, we refer
to the preliminary version \cite{NormalForm} of this paper, section 7.

\begin{proof}[Proof of Corollaries \ref{C:vectfields} and \ref{C:vectfieldsbis}]
We go back to preliminary form
$$X=w\partial_z+(g_0(z)+g_1(z)w)\partial_w$$
(see proof of Theorem \ref{T:vectfields}). Following \cite{Mez} 
(see also \cite{NormalForm}), the foliation
$\mathcal F$ either admits an invariant curve of the form $C:\{w^2+a(x)w+b(x)=0\}$,
where $a(z)$ and $b(z)$ are (real or complex) analytic functions vanishing at $0$,
or admits a smooth (real or complex) analytic invariant curve transversal to the
fibration $\{w=\text{constant}\}$.
We want to simplify this invariant curve by a change of coordinates of the form
$(z,w):=(\varphi(z),w+\phi(z))$. The vector field will therefore
take the more general form 
$$X=(f_0(z)+f_1(z)w)\partial_z+(g_0(z)+g_1(z)w)\partial_w.$$
In the former case, the invariant curve is a $2$-fold covering of the $z$-variable.
One can use a vertical translation
$w:=w+\phi(z)$ so that $C$ becomes invariant by the involution $(z,w)\mapsto(z,-w)$,
i.e. $C=\{w^2=\tilde b(z)\}$. Then, by a change of the $z$-coordinate, one can normalize
$\tilde b(z)=z^k$ (or $\tilde b(z)=-z^k$ when $k$ is even in the real setting).
In these new coordinates, writing that each of the vector fields $X\pm i_*X$ vanish identically along
the curve $t\mapsto(t^k,t^2)$, we deduce that $X$ is takes the form (1)
(or (3) when $k$ is even in the real setting) of Corollary \ref{C:vectfieldsbis}.

\eject

In the saddle case, we have $k=2$.
We set $f(z):={g_0(z)\over g_0(0)}$
and $g(z):=g_1(z)-{g_1(0)\over g_0(0)}g_0(z)$
so that $f(0)=1$, $g(0)=0$ and
the vector field $X$ writes
$$X=f(z)()+g(z)(z\partial_z+w\partial_w)\ \ \
\text{with}\ \ \
X_1=\left(\begin{array}{cc}
  g_1(0) &  g_0(0) \\
  g_0(0) &  g_1(0) \\
\end{array}\right)$$
($g_0(0)=\pm(\lambda_2-\lambda_1)\not=0$).
Finally, after a rotation $(z,w):=(z-w,z+w)$, we obtain normal forms (1)
of Corollary \ref{C:vectfields} for saddles.

In the case $\mathcal F$ admits a smooth analytic invariant curve transversal to the
fibration $\{w=\text{constant}\}$, we first use a vertical translation $w:=w+\phi(z)$
to straighten it
onto the horizontal axis and then use change of $z$-coordinate to send
the tangency set $\{\Sigma\}$ between the foliation $\mathcal F$
and the vertical fibration onto the line $\{w=z\}$.
We immediately obtain normal form (2) of Corollary \ref{C:vectfieldsbis}
 (resp. of Corollary \ref{C:vectfields} in the saddle-node case $k=1$).
\end{proof}

\end{document}